\documentclass[12pt]{amsart}

\usepackage{amssymb}
\usepackage{amscd}
\usepackage[mathscr]{eucal}
\usepackage{amssymb}
\usepackage{xspace}

\theoremstyle{plain}
\newtheorem{theorem}{Theorem}[section]
\newtheorem{corollary}[theorem]{Corollary}
\newtheorem{lemma}[theorem]{Lemma}

\newtheorem{remark}[theorem]{Remark}

\numberwithin{equation}{section}
 
\newcommand{\Q}{\ensuremath{\mathcal{Q}}\xspace}

\newcommand{\Ss}{\ensuremath{\mathcal{S}}\xspace}
\newcommand{\T}{\ensuremath{\mathcal{T}}\xspace}
\newcommand{\K}{\ensuremath{\mathcal{K}}\xspace}

\newcommand{\Rd}{\ensuremath{\mathbb{R}^{d}}\xspace}

\newcommand{\R}{\ensuremath{\mathbb{R}}\xspace}
\newcommand{\N}{\ensuremath{\mathbb{N}}\xspace}
\newcommand{\Z}{\ensuremath{\mathbb{Z}}\xspace}
\newcommand{\Zd}{\ensuremath{\mathbb{Z}^{d}}\xspace}

\begin{document}
\title[The Unit Cube]
   {Spectral and Tiling properties of the Unit Cube$^{\dagger}$}

\author{Alex Iosevich}\author{Steen Pedersen}
\address{Department of Mathematics, Wright State University,
   Dayton OH 45435, USA}
\email{iosevich@zara.math.wright.edu}
\email{steen@math.wright.edu}
\thanks{$^{\dagger}$Research partially
supported by the NSF grant DMS97-06825}

\begin{abstract}
Let $\Q=[0,1)^d$ denote the unit cube in $d$-dimensional Euclidean
space \Rd and let \T be a discrete subset of \Rd.
We show that the exponentials $e_t(x):=exp(i2\pi tx)$, $t\in\T$ form
an othonormal basis for $L^2(\Q)$ if and only if the translates
$\Q+t$, $t\in\T$ form a tiling of \Rd.
\end{abstract}

\maketitle


\section{Introduction}
Let $\Q:=[0,1)^d$ denote the unit cube in $d$-dimensional Euclidean
space \Rd.
Let \T be a discrete subset of \Rd.
We say \T is a \emph{tiling set} for \Q, if each $x\in\Rd$ can be
written uniquely as $x=q+t$, with $q\in\Q$ and $t\in \T$.
We say \T is a \emph{spectrum} for \Q, if the exponentials
\[
   e_t(x):=e^{i2\pi tx}, t\in \T
\]
form an orthonormal basis for $L^2(\Q)$. Here juxtaposition $tx$ of
vectors $t,x$ in \Rd denote the usual inner product
$tx=t_1x_1+\cdots+t_dx_d$ in \Rd and $L^2(\Q)$ is equipped with the
usual inner product, viz.,
\begin{equation*}
   \langle f, g \rangle := \int_{\Q}f\,\overline{g}\,dm
\end{equation*}
where $m$ denotes Lebesgue measure.
The main result proved in this paper is
\begin{theorem}\label{T:main}
   Let \T be a subset of \Rd. Then \T is a spectrum for the unit cube
   \Q if and only if \T is a tiling set for the unit cube \Q.
\end{theorem}

\begin{remark}
As we shall discuss below there exists highly counter-intuitive
cube-tilings in $\mathbb{R}^d$ for sufficiently large $d$.
Those tilings can be much more complicated than lattice tilings.
Theorem \ref{T:main} is clear if \T is a lattice. The point of Theorem
\ref{T:main} is that the result still holds even if the restrictive
lattice assumption is dropped.
\end{remark}

Sets whose translates tile \Rd and the corresponding tiling sets
have been investigated intensively, see \cite{GN}, \cite{LaWa1},
\cite{LaWa2} for some recent papers.
Even the one-dimensional case $d=1$ is
non-trivial. The study of sets whose $L^2$-space admits orthogonal
bases of exponentials was begun in \cite{Fu}. Several papers have
appeared recently, see e.g.,
\cite{JP1}, \cite{Pede}, \cite{LaWa3}.
It was conjectured in \cite{Fu} that a set admits a tiling set if
and only if it admits a spectrum, i.e., the corresponding $L^2$-space
admits an orthogonal basis of exponentials.

Cube tilings have a long history beginning with a conjecture due to
Minkowski: in every lattice tiling of \Rd by translates of \Q some
cubes must share a complete $(d-1)$-dimensional face.
Minkowski's conjecture was proved in \cite{Haj1}, see \cite{SS} for a
recent exposition. Keller \cite{Kel1}
while working on Minkowski's conjecture made the stronger conjecture
that one could omit the lattice assumption in Minkowski's conjecture.
Using \cite{Sza} and \cite{CoSz} it was shown in \cite{LaSh}
that there are
cube tilings in dimensions $d\geq 10$ not satisfying Keller's
conjecture.

The study of the possible spectra for the unit cube was initiated in
\cite{JP2}, where Theorem \ref{T:main} was conjectured.
Theorem \ref{T:main} was proved in \cite{JP2} if $d\leq 3$
and for any $d$ if \T is periodic.

The terminology \emph{spectrum for \Q} originates in a problem about
the existence of certain commuting self-adjoint partial differential
operators.
We say that two self-adjoint operators \emph{commute} if their
spectral measures commute, see \cite{RS} for an introduction to
the theory of unbounded self-adjoint operators. The following result
was proved in \cite{Fu} under a mild regularity condition on the
boundary, the regularity condition was removed in \cite{Ped1}.
\begin{theorem}
   Let $\Omega$ be a connected open subset of \Rd with
   finite Lebesgue measure.
   There exists a set \T so that the exponentials $e_t$, $t\in\T$
   form an orthogonal basis for $L^2(\Omega)$ if and only if there
   exists commuting self-adjoint operators
   $H=(H_1,\ldots,H_d)$ so that each
   $H_j$ is defined on $C^{\infty}_c(\Omega)$ and
   \begin{equation}\label{E:partial}
      H_j f= \frac{1}{i2\pi}\frac{\partial f}{\partial x_j}
   \end{equation}
   for any $f\in C^{\infty}_c(\Omega)$ and any $j=1,\ldots,d$.
\end{theorem}
More precisely, if $e_t$, $t\in\T$ is an orthogonal basis for
$L^2(\Omega)$ then a commuting tuple $H=(H_1,\ldots,H_d)$ of
self-adjoint operators satisfying (\ref{E:partial}) is uniquely
determined by $H_j e_t = t_j e_t$, $t\in\T$. Conversely, if
$H=(H_1,\ldots,H_d)$ is a commuting tuple of self-adjoint operators
satisfying (\ref{E:partial}) then the joint spectrum $\sigma(H)$ is
discrete and each $t\in\sigma(H)$ is a simple
eigen-value corresponding
to the eigen-vector $e_t$, in particular, $e_t$, $t\in\sigma(H)$
is an orthogonal basis for $L^2(\Omega)$.

We prove any tiling set is a spectrum
in Section \ref{S:tiling}, the
converse is proved in Section \ref{S:s_to_t}.
Key ideas in both proofs are that if $(g_n)$ is an orthonormal
family in
$L^2(\Q)$ and $f\in L^2(\Q)$ then we have equality in Bessel's
inequality
\[
   \sum |\langle f,g_n \rangle |^2 \leq ||f||^2
\]
if and only if $f$ is in the closed linear span of $(g_n)$, and
a sliding lemma (Lemma \ref{L:sliding}) showing that we may
translate certain parts of a spectrum or tiling set while preserving
the spectral respectively the tiling set property.
In Section
\ref{S:spectrum} we prove some elementary properties of spectra and
tiling sets.
For $t\in\Rd$ let $\Q+t:=\{q+t:q\in\Q\}$ denote the translate of \Q
by the vector $t$.
We say $(\Q,\T)$ is \emph{non-overlapping} if the cubes $\Q+t$
and $\Q+t^{\prime}$ are disjoint for any $t,t^{\prime}\in\T$.
Note, \T is a tiling set for \Q if and only if $(\Q,\T)$ is
non-overlapping and $\Rd=\Q_{\T}:=\bigcup_{t\in\T}(\Q+t)$.
We say $(\Q,\T)$ is \emph{orthogonal}, if the exponentials
$e_t, t\in \T$ are orthogonal in $L^2(\Q)$. A set \T is a spectrum
for \Q if and only if $(\Q,\T)$ is orthogonal and
\[
   \sum_{t\in\T} |\langle e_n,e_t \rangle |^2 =1
\]
for all $n\in\Zd$.
Let
$\N$ denote the positive integers $\{1,2,3,\ldots\}$ and let
$\Z$ denote the set of all integers $\{\ldots,-1,0,1,2,\ldots\}$.

As this paper was in the final stages of preparation, we received a
preprint \cite{LRW} by Lagarias, Reed and Wang proving our main
result.Compared to \cite{LRW} our proof that any spectrum is a tiling
set uses completely different techniques, the proof that any tiling
set is a spectrum is similar to the proof in \cite{LRW} in that both
proofs makes use of Keller's Theorem (Theorem \ref{T:keller}) and an
argument involving an inequality becoming equality. We wish to thank
Lagarias for the preprint and useful remarks. Robert S. Strichartz
helped us clarify the exposition.

\section{Plan}\label{S:plan}
Our plan is as follows. The basis property is
equivalent to the statement that the sum
\begin{equation}\label{E:plancerel}
   \sum_{t\in\T} |\langle e_x,e_t \rangle |^2 =1
\end{equation}
for all $x\in\Rd$. It is easy to see that if \T has the basis
property then the cubes $\Q+t$, $t\in\T$ are non-overlapping.
We show by a \emph{geometric argument}
that if the basis
property holds and the tiling
property does not hold then the sum in (\ref{E:plancerel})
is strictly less than one.
Conversely, if \T has the tiling property then
the exponentials $e_t$, $t\in\T$
are orthogonal by Keller's Theorem, Plancerel's
Theorem now implies that the sum in (\ref{E:plancerel}) is one.
The geometric argument is based Lemma \ref{L:sliding}, an analoguous
lemma was used by Perron \cite{Per1} in his proof of Keller's Theorem.


\section{Spectral Properties}\label{S:spectrum}
We begin by proving a simple result characterizing orthogonal subsets
of \Rd. There is a corresponding (non-trivial) result for tilings,
stated as Theorem \ref{T:keller} below.
\begin{lemma}[Spectral version of Keller's theorem]
\label{L:orthogonal}
   Let \T be a discrete subset of \Rd. The pair
   $(\Q,\T)$ is orthogonal if and only if
   given any pair
   $t,t^{\prime}\in\T$, with $t\neq t^{\prime}$, there exists a
   $j\in\{1,\ldots,d\}$ so that $|t_j-t^{\prime}_j|\in\N$.
\end{lemma}
\begin{proof} For $t,t^{\prime}\in\Rd$ we have
   \begin{equation}
      \langle e_t, e_{t^{\prime}} \rangle
      =\prod_{j=1}^d \phi(t_j-t^{\prime}_j)
   \end{equation}
   where for $x\in\R$
   \begin{equation}\label{E:phi}
      \phi(x):=
      \begin{cases}
         1,    &\text{if $x=0$;}\\
         \frac{e^{i2\pi x}-1}{i2\pi x},    &\text{if $x\neq 0$.}
      \end{cases}
   \end{equation}
   The lemma is now immediate.
\end{proof}
We can now state our first result showing that there is a connection
between spectra and tiling sets for the unit cube.
\begin{corollary}\label{C:packing}
   Let \T be a subset of \Rd.
   If $(\Q,\T)$ is orthogonal, then $(\Q,\T)$ is non-overlapping.
\end{corollary}
A key technical lemma needed for our proofs of both implications in
our main result is the following lemma. The lemma shows that
a certain part of
a spectrum (respectively tiling set) can be translated
independently of its
complement without destroying the spectral (respectively tiling set)
property. The tiling set part of the lemma if taken from \cite{Per1}.
\begin{lemma}\label{L:sliding}
   Let \T be a discrete subset of \Rd, fix $a,b\in\R$. Let
   $c:=(b,0,\ldots,0)\in\Rd$ and for $t\in\T$ let
   \begin{equation*}
      \alpha_{\T,a,b}(t)
         :=\begin{cases}
              t, &\text{if $t_1-a\in\Z$}; \\
              t+c, &\text{if $t_1-a\notin\Z$}.
           \end{cases}
   \end{equation*}
   We have the following conclusions:
   (a) If \T is a spectrum for \Q, so is $\alpha_{\T,a,b}(\T)$.
   (b) If \T is a tiling set for \Q, so is $\alpha_{\T,a,b}(\T)$.
\end{lemma}
\begin{proof}
   Suppose \T  is a spectrum for \Q.
   The orthogonality of $(\Q,\alpha_{\T,a,b}(\T))$
   is an easy consequence of Lemma \ref{L:orthogonal}.
   Let $A_{\T,a,b}e_t:=e_{\alpha_{\T,a,b}(t)}$ for $t\in\T$. To
   simplify the notation we will write $A_b$ in place of
   $A_{\T,a,b}$.
   By orthogonality and linearity $A_b$ extends to
   an isometry mapping $L^2(\Q)$ into itself. We must show that
   the range $A_bL^2(\Q)$ is all of $L^2(\Q)$. Let $\K_+$ be the
   subspace of $L^2(\Q)$ spanned by the exponentials $e_t$, $t\in\T$
   with $t_1-a\in\Z$ and let $\K_-$ be the
   subspace of $L^2(\Q)$ spanned by the exponentials $e_t$, $t\in\T$
   with $t_1-a\notin\Z$. Then $A_bf=f$ for all $f\in\K_+$, so
   $A_b\K_+=\K_+$. Since $A_b$ preserves orthogonality,
   $A_b\K_-\subseteq\K_-$. We must show
   $A_b\K_-=\K_-$. Since $b\in\R$ is arbitrary, we also have that
   the map $A_{-b}$ is an isometry mapping $\K_-$ into itself.
   By construction
   $A_bf=e_c\,f$ and $A_{-b}f=\overline{e_c}\,f$
   for all $f\in\K_-$. It follows that
   $\K_-=A_{b}A_{-b}\K_-\subseteq A_{b}\K_-\subseteq\K_-$. Hence,
   $A_b\K_-=\K_-$ as desired. The proof that $\alpha_{\T,a,b}(\T)$ is
   a tiling set provided \T is, follows from the last part of the
   proof of Theorem \ref{T:keller} below.
\end{proof}

\section{Any spectrum is a tiling set}\label{S:s_to_t}
For $n^{\prime}\in\Z^{d-1}$ let $\ell_{n^{\prime}}$ be the line
in \Rd given by $\{(x,n^{\prime})$, $x\in\R\}$.
The idea of our proof that any spectrum for
\Q must be a tiling set for \Q is as follows. Suppose \T is a
spectrum but not a tiling set. Fix $n^{\prime}\in\Z^{d-1}$ and
pick a $t\in\T$ (if any) so that $\Q+t$ intersects the line
$\ell_{n^{\prime}}$ applying Lemma \ref{L:sliding} we
can insure that $t_1\in\Z$. Repeating this for
each $n^{\prime}\in\Z^{d-1}$ we can ensure $t_1\in\Z$ for any
$t\in\T^{\text{new}}$.
Considering each of
the remaining coordinate directions we end up with $\T^{\text{new}}$
being a subset of $\Zd$. (The meaning of $\T^{\text{new}}$ changes
with each application of Lemma \ref{L:sliding}.)
By Lemma \ref{L:contradiction} $\T^{\text{new}}$ is not a tiling
set for
\Q since \T was not a tiling set, so $\T^{\text{new}}$ is a proper
subset of \Zd, contradicting the basis property. The difficulty with
this outline is that after we apply Lemma
\ref{L:sliding} an infinite number of times the basis property may
not hold. In fact, associated to each application of Lemma
\ref{L:sliding} is an isometric isomorphism $A_{b_n}$. Without
restrictions on the sequence $(b_n)$ the
infinite product $\prod_{n=1}^{\infty}A_{b_n}$
need not be convergent
(e.g., with respect to the weak operator topology). Even if the
infinite product $\prod_{n=1}^{\infty}A_{b_n}$ is convergent, the
limit may be a non-surjective isometry.

It turns out that if we use Lemma \ref{L:sliding} to put a
large finite part of \T into \Zd then we can use decay properties of
the Fourier transform of the characteristic function of the cube \Q
to contradict (\ref{E:plancerel}).

The following lemma shows that sums of the Fourier transform of the
characteristic function of the cube \Q over certain discrete sets
has uniform decay properties.
\begin{lemma}\label{L:estimate}
   Let $\phi$ be given by (\ref{E:phi}).
   There exists a constant $C>0$ so that
   \begin{equation*}
      \sum_{t\in\T_N} \prod_{j=1}^d|\phi(t_j)|^2\leq\frac{C}{N}
   \end{equation*}
   for any $N>1$,
   whenever $\T\subset\Rd$ is a spectrum for the unit cube \Q.
   Here $\T_N$ is the set of $t\in\T$ for which $|t_j|>N$,
   for at least one $j$. Note, the constant $C$ is uniform over all
   spectra \T for the unit cube \Q and all $N>1$.
\end{lemma}
\begin{proof}
   Let \T be a spectrum for \Q.
   For any partition $P=\{I,II,III,IV\}$ of $\{1,\ldots,d\}$, let
   $\T_{N,P}$ denote the set of $t\in\T_N$ for so that $t_j>N$ for
   $j\in I$; $t_j<-N$ for $j\in II$; $0\leq t_j\leq N$ for $j\in III$
   and $-N\leq t_j <0$ for $j\in IV$.
   Note $\T_{N,P}$ is empty unless $I\cup II$ is non-empty.
   For $x\in\R$ let $\psi(x)=1$, if
   $-1<x<1$ and let $\psi(x)=x^{-2}$ if $|x|\geq1$. Then for
   $t\in\T_{N,P}$,
   \begin{equation}\label{E:pest}
      \prod_{j=1}^{d}|\phi(t_j)|^2 \leq \prod_{j=1}^{d} \psi(s_j)
   \end{equation}
   for any $s=(s_1,\ldots,s_d)$ in the cube $X_{t,P}$ given by
   $t_j-1\leq s_j < t_j$ if
   $j\in I\cup III$, and $t_j\leq s_j < t_j+1$ if $j\in II\cup IV$.
   It follows from (\ref{E:pest}) and disjointness (Lemma
   \ref{L:orthogonal}) of the cubes $X_{t,P}$, $t\in\T_{N,P}$ that
   \begin{align*}
      \sum_{t\in\T_{N,P}} \prod_{j=1}^d |\phi(t_j)|^2
      \leq
      \sum_{t\in\T_{N,P}}\int_{X_{t,P}}
         \prod_{j=1}^d \psi(s_j)\,ds
      \leq
      \int_{Y_{N,P}}\prod_{j=1}^d \psi(s_j)\,ds,
   \end{align*}
   where $Y_{N,P}$ is the set of $y\in\Rd$ for which $N-1<y_j$ for
   $j\in I$, $y_j<-N+1$ for $j\in II$, $-1<y_j<N$ for $j\in III$, and
   $-N<y_j<1$ for $j\in IV$. By definition of $\psi$ we have
   \begin{equation*}
      \int_{Y_{N,P}}\prod_{j=1}^d \psi(s_j)\,ds
      \leq 3^{d-n}\frac{1}{(N-1)^{n}},
   \end{equation*}
   where $n>0$ is the cardinality of $I\cup II$. Since the number of
   possible partitions $P=\{I,II,III,IV\}$ only depends on the
   dimension $d$ of \Rd, the proof is complete
\end{proof}

The following lemma shows that if \T is a spectrum but not a tiling
set for \Q then the set constructed in Lemma \ref{L:sliding} is also
not a tiling set for \Q. It is needed because the inverse of the
transformation in Lemma \ref{L:sliding} is not of the same form.

\begin{lemma}\label{L:contradiction}
    If \T is a spetrum for \Q but not a tiling set for \Q, then
       $\alpha_{\T,a,b}(\T)$ is not a tiling set for \Q.
\end{lemma}
\begin{proof}
 Suppose \T is a spectrum
   for \Q but not a tiling set for \Q. Let $g\notin\Q_{\T}$. Let
   $\ell:=\{(x,g_2,g_3,\ldots,g_d\}$. If $r,s\in\T$ are so that
   $\Q+r$ and $\Q+s$ intersect $\ell$ then it follows from Lemma
   \ref{L:orthogonal} that $s_1-r_1$ is an integer, since
   $|s_{j}-t_{j}|<1$ for $j\neq 1$ because $\Q+r$ and $\Q+s$
   intersect $\ell$. So either
   $t_1-a\in\Z$ for all $t\in\T$ so that $\Q+t$ intersects $\ell$ or
   $t_1-a\notin\Z$ for all $t\in\T$ so that $\Q+t$ intersects $\ell$.
   In the first case $g\notin\Q_{\alpha_{\T,a,b}(\T)}$ in the second
   case $g+c\notin\Q_{\alpha_{\T,a,b}(\T)}$.
\end{proof}

\begin{proof}[Proof of basis implies tiling]
   Suppose \T is a spectrum for the unit cube \Q. By Corollary
   \ref{C:packing} the pair $(\Q,\T)$ is non-overlapping. We must
   show that the union $\Q_{\T}=\cup_{t\in\T}(\Q+t)$ is all of \Rd.
   To get
   a contradiction suppose $g\notin\Q_{\T}$.
   Let $N$ be so large that
   $g\in(-N+2,N-2)^d$. Let $\T(N):=\T\cap(-N-1,N+1)^d$.
 
   Let $n^{\prime}_1:=(-N,-N,\ldots,-N)\in\Z^{d-1}$.
   Pick $t\in\T(N)$ so that $\Q+t$ intersects $\ell_{n^{\prime}_1}$
   (if such a $t$ exists). Use Lemma \ref{L:sliding} with $a=0$ and
   $b=b_1:=t_1-\lfloor t_1 \rfloor$ to conclude
   $\T_1:=\alpha_{\T,a,b}(\T)$ has the basis property. It follows
   from Lemma \ref{L:orthogonal} that
   $t_1\in\Z$ for any $t\in\T_1$ so that
   $\Q+t$ intersects $\ell_{n^{\prime}_1}$.
 
   Let $n^{\prime}_2:=(-N,-N,\ldots,-N,-N+1)\in\Z^{d-1}$.
   Pick $t\in\T_{1}(N)$ so that $\Q+t$ intersects
   $\ell_{n^{\prime}_2}$ (if such a $t$ exists).
   Use Lemma \ref{L:sliding} with $a=0$ and
   $b=b_2:=t_1-\lfloor t_1 \rfloor$ if
   $b_1+t_1-\lfloor t_1\rfloor\geq 1$
   and $b=b_2:=t_1-\lfloor t_1 \rfloor-1$ if
   $b_1+t_1-\lfloor t_1\rfloor<1$ to conclude
   $\T_2:=\alpha_{\T_1,a,b}(\T_1)$ has the basis property.
   It follows from Lemma \ref{L:orthogonal} that
   $t_1\in\Z$ for any $t\in\T_2$ so that $\Q+t$ intersects
   $\ell_{n^{\prime}_2}$. Note we did not move any of the cubes in
   $\T_1$ with $-N-1<t_j\leq-N$, for $j=2,\ldots,d$.

   Continuing in this manner, we end up with $\T^{\prime}$ having
   the basis property so that $t_1\in\Z$ for any $t\in\T^{\prime}$
   with $-N-1<t_j<N+1$ for $j=2,\ldots,d$. Note $-1<\sum_1^n b_j<1$
   for any $n$. So if at some stage $t\in\T_n$ is derived from
   $t^{\text{original}}\in\T$ then we have
   $t^{\text{original}}_1-1<t_1<t^{\text{original}}_1+1$.
 
   Repeating this process for each of the other coordinate directions
   we end up with $\T^{\text{new}}$ so that $\T^{\text{new}}(N-1)$ is
   a subset of the integer lattice \Zd,
   any $t\in\T^{\text{new}}(N-1)$ is
   obtained from some $t^{\text{original}}\in\T(N)$, and any
   $t^{\text{original}}\in\T(N-1)$ is translated onto some
   $\T^{\text{new}}(N)$. In short, we did not move any point in \T
   very much. By Lemma \ref{L:contradiction} it follows that
   $(-N,N)^d\setminus\Q_{\T^{\text{new}}}$ is non-empty, hence there
   exists $g^{\text{new}}\in\Zd$, so that
   $g^{\text{new}}\in(-N,N)^d\setminus\T^{\text{new}}$. Replacing
   $\T^{\text{new}}$ by $\T^{\text{new}}-g^{\text{new}}$,
   if necessary, and
   applying the process described above we may assume
   $g^{\text{new}}=0$. To simplify the notation let
   $\T=\T^{\text{new}}$. We have
   \begin{align*}
      1=\sum_{t\in\T} |\langle e_t,e_0 \rangle|^2
       =\sum_{t\in\T(N)} |\langle e_t,e_0 \rangle|^2
          +\sum_{t\in\T_N} |\langle e_t,e_0 \rangle|^2.
   \end{align*}
   The first sum $=0$ since $\T(N)\subset\Zd$ and $0\notin\T(N)$, the
   second sum is $<1$ for $N$ sufficiently large
   by Lemma \ref{L:estimate}.
   This contradiction completes the proof.
\end{proof}


\section{Any tiling set is a spectrum}\label{S:tiling}
The following result (due to \cite{Kel1})
shows that any tiling set for the cube is
orthogonal. It is a key step in our proof that any tiling set for the
cube must be a spectrum for the cube and should be compared with
Lemma \ref{L:orthogonal} above.
The proof is essentially taken from \cite{Per1}.
\begin{theorem}[Keller's Theorem]\label{T:keller}
   If \T is a tiling set for \Q, then given any pair
   $t,t^{\prime}\in\T$, with $t\neq t^{\prime}$, there exists a
   $j\in\{1,\ldots,d\}$ so that $|t_j-t^{\prime}_j|\in\N$.
\end{theorem}
\begin{proof}
   Let \T be a tiling set for \Q. Suppose $t,t^{\prime}\in\T$.
   The proof is by induction on the number of $j$'s for which
   $|t_{j}-t_{j}^{\prime}|\geq 1$. Suppose
   that $|t_j-t_j^{\prime}|<1$ for all but one $j\in\{1,\ldots,d\}$.
   Let $j_0$ be the exceptional $j$, then
   $|t_{j_0}-t_{j_0}^{\prime}|\geq 1$. Fix $x_j$, $j\neq j_0$ so
   that the line $\ell_{j_0}:=\{(x_1,\ldots,x_d):x_{j_0}\in\R\}$
   passes through both of the cubes $\Q+t$ and $\Q+t^{\prime}$.
   Considering the cubes $Q+t$, $t\in\T$ that intersect $\ell_{j_0}$
   it is immediate that $|t_{j_0}-t_{j_0}^{\prime}|\in\N$.
 
   For the inductive step,
   suppose $|t_j-t_j^{\prime}|<1$ for $k$ values of $j$ and
   $|t_j-t_j^{\prime}|\geq 1$ for the remaining $d-k$ values of $j$
   implies $|t_{j_0}-t_{j_0}^{\prime}|\in\N$ for some $j_0$. Let
   $t,t^{\prime}\in\T$ be so that $|t_j-t_j^{\prime}|<1$ for $k-1$
   values of $j$ and $|t_j-t_j^{\prime}|\geq 1$ for the remaining
   $d-k+1$ values of $j$. Interchanging the coordinate axes, if
   necessary, we may assume
   \begin{align*}
      |t_j-t_j^{\prime}|&\geq 1,\text{ for } j=1,\ldots,d-k+1 \\
      |t_j-t_j^{\prime}|&<1, \text{ for } j=d-k+2,\ldots,d.
   \end{align*}
   If $t_1-t_1^{\prime}$ is an integer, then there we are done.
   Assume $t_1-t_1^{\prime}\notin\Z$. Let
   $c:=(t_1-t_1^{\prime},0,\ldots,0)$, and
   for $\tilde{t}\in\T$ let
   \begin{equation*}
      s(\tilde{t}):=\begin{cases}
         \tilde{t}-c,
            &\text{ if $\tilde{t}_1-t_1\in\Z$} \\
         \tilde{t},
            &\text{ if $\tilde{t}_1-t_1\notin\Z$}.
            \end{cases}
   \end{equation*}
   In particular, $s(t)=t-c$ and
   $s(t^{\prime})=t^{\prime}$.
   We claim the set $\Ss:=\{s(\tilde{t}):\tilde{t}\in\T\}$ is a
   tiling set for \Q. Assuming, for a moment, that the claim is
   valid,
   we can easily complete the proof.
   In fact,
   $|s(t)_1-s(t^{\prime})_1|=0$ and $|s(t)_j-s(t^{\prime})_j|<1$
   for
   $j=d-k+2,\ldots,d$, so by the inductive hypothesis one of the
   numbers $t_j-t_j^{\prime}=s(t)_j-s(t^{\prime})_j$,
   $j=2,\ldots,d-k+1$ is a non-zero integer.
 
   It remains to prove that $\Ss$ is a tiling set for \Q. We must
   show that $\Q_{\Ss}$ is non-overlapping and that
   $\Rd\subset\Q_{\Ss}$. First we dispense with the non-overlapping
   part. Let $a,a^{\prime}$ be distinct points in \T.
   Suppose $x$ is a point in the intersection
   $(\Q+s(a))\cap(\Q+s(a^{\prime}))$, then
   $x-s(a),x-s(a^{\prime})\in\Q$, in particular, $0\leq x_j-a_j<1$
   and $0\leq x_j-a^{\prime}_j<1$ for $j=2,\ldots,d$. It follows that
   $|a_j-a^{\prime}_j|<1$ for $j=2,\ldots,d$, so first paragraph of
   the proof shows that $|a_1-a^{\prime}_1|\in\N$, hence either
   $a_1-t_1,a_{1}^{\prime}-t_1\in\Z$ or
   $a_1-t_1,a_{1}^{\prime}-t_1\notin\Z$. In both cases we get a
   contradiction to the non-overlapping property of $\Q_{\T}$. In
   fact, if $a_1-t_1,a^{\prime}-t_1\in\Z$, then
   $(\Q+s(a))\cap(\Q+s(a^{\prime}))
   =((\Q+a)\cap(\Q+a^{\prime}))-c=\emptyset$.
   If $a_1-t_1,a_1^{\prime}-t_1\notin\Z$, then
   $(\Q+s(a))\cap(\Q+s(a^{\prime}))
   =((\Q+a)\cap(\Q+a^{\prime}))=\emptyset$.
 
   Let $x\in\Rd$ be an arbitrary point, then $x\in\Q_{\T}$. If
   $x\in (\Q+a)$ for some $a\in\T$ with $a_1-t_1\notin\Z$ then there
   is nothing to prove. Assume $x\in (\Q+a)$ for some $a\in\T$ with
   $a_1-t_1\in\Z$. The point $x+c$ is in $\Q+b$ for some $b\in\T$.
   First we show that $|a_1-b_1|\in\N$. Since $x\in\Q+a$ and
   $x+c\in\Q+b$ we have
   \begin{equation}\label{E:inequalities}
      0\leq x_j-a_j<1, \quad
      0\leq x_j-b_j+c_j<1,
   \end{equation}
   for $j=1,\ldots,d$; so using $c_j=0$, for $j=2,\ldots,d$, it
   follows that $|a_j-b_j|<1$, for $j=2,\ldots,d$; an application
   for the first paragraph of the proof yields the desired result
   that $|a_1-b_1|\in\N$. Using $a_1-t_1\in\Z$ we conclude
   $b_1-t_1\in\Z$; so using the second half of (\ref{E:inequalities})
   and the definition
   of $s(b)$ we have $x\in\Q+s(b)$ as needed.
\end{proof}

\begin{corollary}
   If \T is a tiling set for \Q, then $(\Q,\T)$ is orthogonal.
\end{corollary}
\begin{proof}
   This is a direct consequence of Keller's Theorem and Lemma
   \ref{L:orthogonal}.
\end{proof}
It is now easy to complete the proof that any tiling set for the unit
cube \Q must be a spectrum for \Q.
\begin{proof}[Proof of tiling implies basis]
   Suppose \T is a tiling set for \Q. By Keller's Theorem
   $\{e_t:t\in\T\}$ is an orthogonal set of unit vectors in
   $L^2(\Q)$, so by
   Bessel's inequality
   \begin{equation}\label{E:bessel}
      \sum_{t\in\T} |\langle e_s,e_t\rangle|^2 \leq 1
   \end{equation}
   for any $s\in\Rd$. Note that $\langle e_s,e_t\rangle$ is the
   Fourier transform of the characteristic function of the cube \Q at
   the point $s-t$.
   For any $r\in\Rd$ we have
   \begin{align*}
      1=\int_{\Rd}|\langle e_y,e_0\rangle|^2\, dy
       =\int_{\Q+r}\sum_{t\in\T}|\langle e_x,e_t\rangle|^2\, dx
       \leq \int_{\Q+r}1\,dy =1,
   \end{align*}
   where we used Plancherel's Theorem, the tiling property,
   and Bessel's inequality (\ref{E:bessel}). It follows that
   \begin{equation}\label{E:almost}
      \sum_{t\in\T} |\langle e_s,e_t\rangle|^2=1
   \end{equation}
   for almost every $s$ in $\Q+r$, and since $r$ is arbitrary,
   for almost every $s$ in \Rd. Hence for almost every $s\in\Rd$ the
   exponential $e_s$ is in the closed span of the
   $e_t$, $t\in\T$. This completes the proof.
\end{proof}


\end{document}